\documentclass[review]{elsarticle}

\usepackage{lineno,hyperref}
\modulolinenumbers[5]

\begin{document}

\begin{frontmatter}

\title{Uniformly most powerful unbiased test for conditional independence in Gaussian graphical model}

\author{Koldanov Petr \footnote{Corresponding author, pkoldanov@hse.ru}, Koldanov Alexander, Kalyagin Valeriy; Pardalos Panos } %
\address{National Research University Higher School of Economics, Laboratory of Algorithms and Technologies for Network Analysis, Rodionova 136, Nizhny Novgorod, Russia; Preeminent Professor
 in Industrial and Systems Engineering, Director, CAO  Industrial and Systems Engineering, University of Florida, Gainesville
}




\begin{abstract}
Model selection for Gaussian concentration graph is based on multiple testing
of pairwise conditional independence. In practical applications partial correlation tests are widely used. However it is not known whether partial correlation test is uniformly most
powerful for pairwise conditional independence testing. This question is answered in the paper. Uniformly most
powerful unbiased test of  Neymann structure is obtained. It turns
out, that this test can be reduced to usual partial correlation test. It implies that
partial correlation test is uniformly most powerful unbiased one.
\end{abstract}

\begin{keyword}
\texttt Conditional independence \sep Exponential families \sep Multivariate normal distribution  \sep Sample partial correlation test \sep Tests of Neyman structure \sep Uniformly most powerful unbiased tests.
\end{keyword}

\end{frontmatter}

\linenumbers

\section{Introduction.}
Let $X=(X_1,X_2,\ldots,X_N)$ be a random vector with multivariate Gaussian distribution.
Concentration graph is defined as follows: nodes of the graph are associated with random variables $X_1,X_2,\ldots,X_N$, edge $(i,j)$ is included in the graph iff random variables $X_i,X_j$ are conditionally dependent \cite{bib1}, \cite{bib7}. Model selection for Gaussian concentration graph consists of identification of concentration graph from observations. This problem has a practical importance  in biology and genetics \cite{bib2}, \cite{bib4}. Common approach for model selection is based on multiple testing of individual hypotheses of pairwise conditional independence \cite{bib5}, \cite{bib6}.
 
Conditional independence of $X_i$, $X_j$ given $X_k$, $k \in N(i,j) = \{1,2,\ldots,N\}\setminus \{i,j\}$ is equivalent to the equation $\rho_{i,j \bullet N(i,j)}=0$, where $\rho_{i,j \bullet N(i,j)}$ is the partial correlation of $X_i$ and $X_j$ given $X_k$, $k \in N(i,j)$. For testing hypothesis $\rho_{i,j \bullet N(i,j)}=0$ for multivariate normal distributions the test of sample partial correlation is largely used \cite{bib7}. At the same time, as far as we know, there are no results concerning uniformly most powerful unbiased (UMPU) tests for conditional independence. Such test is of own interest and could improve multiple testing procedures for model selection. In the present paper we construct a uniformly most powerful unbiased test of  Neyman structure for testing pairwise conditional independence. It turns out that this test can be reduced to the sample partial correlation test. Therefore the sample partial correlation test is uniformly most powerful unbiased one. This fact has some important consequences for multiple testing with additive loss function.

The paper is organized as follows. Section \ref{Basic notations and problem statement} contains basic definition and problem statement. In  Section \ref{Tests of Neyman structure} a general description of the tests of Neyman structure is given. In Section \ref{Uniformly most powerful unbiased test} the UMPU  test for testing pairwise conditional independence is constructed. In Section \ref{Sample partial correlation test} it is proved that the UMPU test can be reduced to the sample partial correlation test. 

\section{Basic notations and problem statement.}\label{Basic notations and problem statement}

Let random vector $X=(X_1,X_2,\ldots,X_N)$ have a multivariate normal distribution $N(\mu, \Sigma)$, where $\mu=(\mu_1,\mu_2,\ldots, \mu_N)$ is the vector of means and 
$\Sigma=(\sigma_{i,j})$ is the covariance matrix, $\sigma_{i,j}=\mbox{cov}(X_i,X_j)$, $i,j=1,2,\ldots,N$.  Let $x(t)$, $t=1,2,\ldots,n$  be a sample of the size $n$ from the distribution of $X$ and
$$ 
s_{i,j}=\frac{1}{n} \Sigma_{t=1}^n(x_{i}(t)-\overline{x_{i}})(x_{j}(t)-\overline{x_{j}}), 
$$
be the sample covariance between $X_i$, $X_j$, where $\overline{x_i}=(1/n)\sum_{t=1}^n x_i(t)$. Denote by $S=(s_{i,j})$ the matrix of sample covariances.  

The inverse matrix for $\Sigma$,  $\Sigma^{-1}=(\sigma^{i,j})$ is known as the concentration or precision matrix for the distribution of $X$.   
For simplicity we use the notation $\rho^{i,j}=\rho_{i,j \bullet N(i,j)}$. The problem of pairwise conditional independence testing has the form:
\begin{equation}\label{main_problem}
h_{i,j}:\rho^{i,j}=0 \mbox{ vs } k_{i,j}:\rho^{i,j}\neq 0
\end{equation}
According to \cite{bib1} the partial correlation can be calculated as  
$$
\rho^{i,j}=-\frac{\sigma^{i,j}}{\sqrt{\sigma^{i,i}\sigma^{j,j}}}
$$
Therefore the problem of pairwise conditional independence testing (\ref{main_problem}) can be formulated as
\begin{equation}\label{Individual hypotheses}
h_{i, j}:\sigma^{i,j}=0, \  \mbox{ vs } \  k_{i, j}:\sigma^{i,j}\neq 0.
\end{equation}

\section{Test of  Neyman structure.}\label{Tests of Neyman structure}

To construct UMPU test for the problem (\ref{Individual hypotheses}) we use a test of Neyman structure for natural parameters of exponential family \cite{bib8}. 
Let $f(x;\theta)$ be the density of the  exponential family:
\begin{equation} \label{exp}
f(x;\theta)=c(\theta)exp\left(\sum_{j=1}^M\theta_jT_j(x)\right)m(x)
\end{equation}
where $c(\theta)$ is a function defined in the parameters space,
$m(x)$, $T_j(x)$ are functions defined in the sample space, and
$T_j(X)$ are the sufficient statistics for $\theta_j,j=1,\ldots,M$.

Suppose that  hypothesis  has the form:
\begin{equation}\label{generating hypotheses}
h_j:\theta_j=\theta^0_j \mbox{  vs  } k_j:\theta_j \neq\theta^0_j, 
\end{equation}
where $\theta^0_j$ is  fixed. 

The UMPU test for hypotheses (\ref{generating hypotheses}) is (see \cite{bib8}, Ch. 4, theorem 4.4.1):
\begin{equation}\label{tsn}
\varphi_j= \left\{
\begin{array}{lc}
0, & \mbox{if} \ \ c'_j(t_1,\ldots,t_{j-1},t_{j+1},\ldots,t_M)<t_j<c''_j(t_1,\ldots,t_{j-1},t_{j+1},\ldots,t_M)\\
1, & \mbox{  otherwise  }
 \end{array}
\right.
\end{equation}
where $t_i=T_i(x),i=1,\ldots,M$. The constants $c'_j$, $c''_j$ are defined from the equations
\begin{equation}\label{neymstruc1}
\int_{c'_j}^{c''_j}f(t_j;\theta^0_j|T_i=t_i,i=1,\ldots,M;i \neq j)dt_j=1-\alpha
\end{equation}
and
\begin{equation}\label{neymstruc2}
\begin{array}{l}
\int_{-\infty}^{c'_j}t_jf(t_j;\theta^0_j|T_i=t_i,i=1,\ldots,M;i\neq j)dt_j+\\
+\int_{c''_j}^{+\infty}t_jf(t_j;\theta^0_j|T_i=t_i,i=1,\ldots,M;i\neq j)dt_j=\\
=\alpha \int_{-\infty}^{+\infty}t_jf(t_j;\theta^0_j|T_i=t_i,i=1,\ldots,M;i\neq
j)dt_j
\end{array}
\end{equation}
where $f(t_j;\theta^0_j|T_i=t_i,i=1,\ldots,M;i\neq j)$ is the
density of conditional distribution of statistic $T_j$  given $T_i=t_i$, $i=1,2,\ldots,N$, $i \neq j$, and $\alpha$ is the significance  level of
the test.

\section{Uniformly most powerful unbiased test for conditional independence.}\label{Uniformly most powerful unbiased test}
Now we construct the UMPU test for testing hypothesis of conditional independence (\ref{Individual hypotheses}). Consider statistics 
$$
S_{k,l}=\frac{1}{n} \Sigma_{t=1}^n(X_{i}(t)-\overline{X_{i}})(X_{j}(t)-\overline{X_{j}}),
$$
Joint distribution  of statistics $S_{k,l}$, $k,l = 1,2,\ldots,N$, $n>N$ is given by Wishart density function \cite{bib7}:
$$
f(\{s_{k,l}\})=\displaystyle \frac{ [\det (\sigma^{k,l})]^{n/2}
\times [\det(s_{k,l})]^{(n-N-2)/2}\times \exp[-(1/2)\sum_k \sum_l
s_{k,l} \sigma^{k,l}]} {2^{(Nn/2)}\times \pi^{N(N-1)/4} \times
\Gamma(n/2)\Gamma((n-1)/2)\cdots\Gamma((n-N+1)/2)}
$$
if the  matrix $S=(s_{k,l})$ is positive definite, and $f(\{s_{k,l}\})=0$ otherwise. It implies that statistics $S_{k,l}$ 
are sufficient statistics for natural parameters $\sigma^{k,l}$. Wishart density function can be written as:
$$
f(\{s_{k,l}\})=\displaystyle C(\{\sigma^{k,l}\}) 
\exp[-\sigma^{i,j}s_{i,j} - \frac{1}{2} \sum_{(k,l)\neq
(i,j);(k,l)\neq(j,i)} s_{k,l} \sigma^{k,l}]  m(\{s_{k,l}\})
$$
where
$$
C(\{\sigma^{k,l}\})=c_1^{-1}[\det (\sigma^{k,l})]^{n/2}
$$
$$
c_1=2^{(Nn/2)}\times \pi^{N(N-1)/4}\times
\Gamma(n/2)\Gamma((n-1)/2)\cdots\Gamma((n-N+1)/2)
$$
$$ m(\{s_{k,l}\})=[\det(s_{k,l})]^{(n-N-2)/2}
$$

According to (\ref{tsn}) the UMPU test for hypothesis (\ref{Individual hypotheses}) has the form:
\begin{equation}\label{Nstructure}
\varphi_{i, j}(\{s_{k, l}\})=\left\{\begin{array}{rl}
 \ 0, &\mbox{}\: if \:  c_{i,j}'(\{s_{k,l}\})<s_{i,j}<c_{i,j}'' (\{s_{k,l}\}),\  (k,l)\neq (i,j)\\
 \ 1, &\mbox{}\: if \: s_{i,j}\leq c_{i,j}'(\{s_{k,l}\})\mbox{ or } s_{i,j}\geq c_{i,j}''(\{s_{k,l}\}),\  (k,l)\neq (i,j)
 \end{array}\right.
\end{equation}
where  the critical values $c'_{i,j}, c''_{i,j}$  are defined from the equations (according to (\ref{neymstruc1}),(\ref{neymstruc2}))
\begin{equation}\label{threshold1}
\displaystyle \frac{\int_{I \cap [c_{i,j}';c_{i,j}'']}
 [\det(s_{k,l})]^{(n-N-2)/2}  ds_{i,j}}
{\int_{I}  [\det(s_{k,l})]^{(n-N-2)/2}
ds_{i,j}} =1-\alpha
\end{equation}
\begin{equation}\label{threshold2}
\begin{array}{l}
\displaystyle \int_{I \cap (-\infty;c_{i,j}']}
s_{i,j}[\det(s_{k,l})]^{(n-N-2)/2}
ds_{i,j}+\\
+\displaystyle \int_{I \cap [c_{i,j}'';+\infty)}
s_{i,j} [\det(s_{k,l})]^{(n-N-2)/2}
ds_{i,j}=\\
 =\alpha \int_I s_{i,j}[\det(s_{k,l})]^{(n-N-2)/2} ds_{i,j}
\end{array}
\end{equation}
where $I$ is the interval of values of $s_{i,j}$ such that the
matrix $S=(s_{k,l})$ is positive definite and  $\alpha$ is the
 significance level of the test.

Let $S=(s_{k,l})$ be positive definite (this is true with probability 1 if $n>N$).  Consider $\det(s_{k, l})$ as a function of the variable $s_{i,j}$ only, when fixing the values of all others $\{s_{k,l}\}$. This determinant is a quadratic polynomial of $s_{i,j}$:
\begin{equation}\label{determinant_equation}
\det(s_{k,l})=-as_{i,j}^2+bs_{i,j}+c
\end{equation} 

Let $K=(n-N-2)/2$. Denote by $x_1,x_2$ ($x_1<x_2$) the roots of the equation $-ax^2+bx+c=0$. One has 
$$
\int_f^d(ax^2-bx-c)^Kdx=\displaystyle (-1)^Ka^K(x_2-x_1)^{2K+1}\int_{\frac{f-x_1}{x_2-x_1}}^{\frac{d-x_1}{x_2-x_1}}u^K(1-u)^Kdu
$$
Therefore the equation (\ref{threshold1}) takes the form:
\begin{equation}\label{threshold_neyman_structure_1}
\displaystyle \int_{\frac{c'-x_1}{x_2-x_1}}^{\frac{c''-x_1}{x_2-x_1}}u^K(1-u)^Kdu=(1-\alpha)\int_0^1 u^K(1-u)^K du
\end{equation}
or
\begin{equation}\label{beta_function}
\displaystyle \frac{\Gamma(2K+2)}{\Gamma(K+1)\Gamma(K+1)}\int_{\frac{c'-x_1}{x_2-x_1}}^{\frac{c''-x_1}{x_2-x_1}}u^K(1-u)^Kdu=(1-\alpha)
\end{equation}
It means that conditional distribution of $S_{i,j}$ when all other $S_{k,l}$ are fixed, $S_{k,l}=s_{k,l}$ is the beta distribution $Be(K+1,K+1)$.

Beta distribution $Be(K+1,K+1)$ is symmetric with respect to the point $\frac{1}{2}$. Therefore the significance level condition (\ref{threshold1}) 
and unbiasedness condition (\ref{threshold2}) are satisfied if and only if:
$$
\displaystyle \frac{c''-x_1}{x_2-x_1}=1 - \frac{c'-x_1}{x_2-x_1}
$$  

Let $q$ be the $\frac{\alpha}{2}$-quantile of beta distribution $Be(K+1,K+1)$, i.e. $F_{Be}(q)=\frac{\alpha}{2}$. Then thresholds  $c'$, $c''$ are defined by: 
\begin{equation}\label{thresholds_from_beta_distr}
\begin{array}{c}
c'=x_1+(x_2-x_1)q \\
c''=x_2-(x_2-x_1)q
\end{array}
\end{equation}

Finally, the UMPU test for testing conditional independence of $X_i$, $X_j$ has the form 
\begin{equation}\label{Neyman_structure_q}
\varphi_{i,j}=\left\{\ 
\begin{array}{ll} 
0, & \displaystyle 2q-1 < \frac{as_{i,j}-\frac{b}{2}}{\sqrt{\frac{b^2}{4}+ac}} < 1-2q \\
1, & \mbox{otherwise}
\end{array}\right.
\end{equation}
where  $a, b, c$ are defined in (\ref{determinant_equation}). 

\section{Sample partial correlation test.}\label{Sample partial correlation test}

It is known \cite{bib1} that hypothesis $\sigma^{i,j}=0$ is equivalent to the hypothesis  $\rho^{i,j}=0$, 
where $\rho^{i,j}$ is the partial correlation between $X_i$ and $X_j$ given $X_k$, $k \in N(i,j)=\{1,2,\ldots,N\}\setminus \{i,j\}$:
$$
\rho^{i,j}=-\frac{\sigma^{i,j}}{\sqrt{\sigma^{i,i}\sigma^{j,j}}}=\frac{- \Sigma^{i,j}}{\sqrt{\Sigma^{i,i}\Sigma^{j,j}}}
$$
where for a given matrix $A=(a_{k,l})$ we denote by $A^{i,j}$  the cofactor of the element $a_{i,j}$. 
Denote by $r^{i,j}$ sample partial correlation 
$$
r^{i,j}=\frac{-S^{i,j}}{\sqrt{S^{i,i}S^{j,j}}}
$$
where $S^{i,j}$ is the cofactor of the element $s_{i,j}$ in the matrix $S$ of sample covariances.

Well known sample partial correlation test for testing hypothesis $\rho^{i,j}=0$ has the form \cite{bib7}:
\begin{equation}\label{Partial_correlation_test}
\varphi_{i,j}=\left\{\ 
\begin{array}{ll} 
0,&|r^{i,j}|\leq c_{i,j}\\
1,&|r^{i,j}|> c_{i,j}
\end{array}\right.
\end{equation} 
where 
$c_{i,j}$ is $(1-\alpha/2)$-quantile of the distribution with the following density function
$$
f(x)=\displaystyle \frac{1}{\sqrt{\pi}}\frac{\Gamma(n-N+1)/2)}{\Gamma((n-N)/2)}(1-x^2)^{(n-N-2)/2}, \ \ \ -1 \leq x \leq 1
$$
Note, that in practical applications  the following Fisher transformation is applied:
$$
z_{i,j}=\frac{\sqrt{n}}{2} \ln\left(\frac{1+r^{i,j}}{1-r^{i,j}}\right)
$$
Under condition $\rho^{i,j}=0$ statistic $Z_{i,j}$  has asymptotically standard normal distribution. 
That is why the following test is largely used in applications \cite{bib4}, \cite{bib5}, \cite{bib6}:
\begin{equation}\label{Fisher_test}
\varphi_{i,j}=\left\{\ 
\begin{array}{ll} 
0,&|z_{i,j}|\leq c_{i,j}\\
1,&|z_{i,j}|> c_{i,j}
\end{array}\right.
\end{equation} 
where the constant $c_{i,j}$ is $(1-\alpha/2)$-quantile  of standard normal distribution. 

In this section we prove that the UMPU test  (\ref{Neyman_structure_q}) can be reduced to the sample partial correlation test (\ref{Partial_correlation_test}), 
and therefore the well known sample partial correlation test for conditional independence is the UMPU one.

\noindent
{\bf Theorem:} {\it Sample partial correlation test (\ref{Partial_correlation_test}) is equivalent to UMPU test (\ref{Neyman_structure_q}) for testing hypothesis 
$\rho^{i,j}=0$ vs $\rho^{i,j} \neq 0$.}

\noindent
{\it Proof:} it is sufficient to prove that
\begin{equation}\label{equality}
\displaystyle \frac{S^{i,j}}{\sqrt{S^{i,i}S^{j,j}}}=\frac{as_{i,j}-\frac{b}{2}}{\sqrt{\frac{b^2}{4}+ac}}
\end{equation}
To prove this equation we introduce some notations. Let $A=(a_{k,l})$ be an $(N \times N)$ symmetric matrix. Fix $i<j$, $i,j =1,2,\ldots,N$.  
Denote by $A(x)$ the matrix obtained from $A$ by replacing the elements $a_{i,j}$ and $a_{j,i}$  by $x$. 
Denote by $A^{i,j}(x)$ the cofactor of the element $(i,j)$ in the matrix $A(x)$. Then the following statement is true

\noindent
{\bf Lemma:} {\it One has $[\mbox{det}A(x)]' = -2A^{i,j}(x)$.}

\noindent
{\it Proof of the Lemma:}  one has from the general Laplace decomposition of $\det A(x)$ by two rows $i$ and $j$:
$$
\det(A(x))=\det \left( 
\begin{array}{ll}
a_{i,i} & x \\
x & a_{j,j} \\
\end{array}  
\right) A^{\{i,j\},\{i,j\}} + \sum_{k<j, k \neq i} \det \left( 
\begin{array}{ll}
a_{i,k} & x \\
a_{j,k} & a_{j,j} \\
\end{array}  
\right)A^{\{i,j\},\{k,j\}}+
$$
$$ + \sum_{k>j} \det \left( 
\begin{array}{ll}
x & a_{i,k} \\
a_{j,j} & a_{j,k} \\
\end{array}  
\right)A^{\{i,j\},\{j,k\}}+ \sum_{k<i} \det \left( 
\begin{array}{ll}
a_{i,k} & a_{i,i} \\
a_{j,k} & x \\
\end{array}  
\right)A^{\{i,j\},\{k,i\}}+ 
$$
$$
\sum_{k>i, k \neq j} \det \left( 
\begin{array}{ll}
a_{i,i} & a_{i,k}\\
x & a_{j,k} \\
\end{array}  
\right)A^{\{i,j\},\{i,k\}}+\sum_{k<l, k,l \neq i,j} \det \left( 
\begin{array}{ll}
a_{i,k} & a_{i,l}\\
a_{j,k} & a_{j,l} \\
\end{array}  
\right)A^{\{i,j\},\{k,l\}}
$$
where $A^{\{i,j\},\{k,l\}}$ is the cofactor of the matrix 
$\left( 
\begin{array}{ll}
a_{i,k} & a_{i,l}\\
a_{j,k} & a_{j,l} \\
\end{array}  
\right)
$
in the matrix $A$. Taking the derivative of $\det A(x)$ one get
$$
[\det(A(x))]'=-2xA^{\{i,j\},\{i,j\}}-\sum_{k<j, k \neq i} a_{j,k}A^{\{i,j\},\{k,j\}}+
\sum_{k>j} a_{j,k}A^{\{i,j\},\{j,k\}}+
$$
$$
+\sum_{k<i} a_{k,i}A^{\{i,j\},\{k,i\}}-
\sum_{k>i, k \neq j} a_{k,i}A^{\{i,j\},\{k,i\}}=-2A^{i,j}(x)
$$
The last equation follows from the symmetry conditions $a_{k,l}=a_{l,k}$ and from Laplace decompositions of $A^{i,j}(x)$ by the row $j$ and the column $i$. 
Lemma is proved. Note, that similar result is proved in (\cite{bib7}, Appendix A).

Now we come back to the proof of the theorem. One has $\det(S(x))=-ax^2+bx+c$, where $a,b,c$ are the same as in (\ref{determinant_equation}). Therefore by Lemma one has $[\det S(x)]'=-2ax+b=-2S^{i,j}(x)$, i.e. 
$S^{i,j}(x)=ax-b/2$. Let 
$x=s_{i,j}$ then $as_{i,j}-\frac{b}{2}=S^{i,j}$.  To prove the theorem it is sufficient to prove that $\sqrt{S^{i,i}S^{j,j}}=\sqrt{\frac{b^2}{4}+ac}$. 
Let $x_2=\frac{b+\sqrt{b^2+4ac}}{2a}$ be the maximum root of equation $ax^2-bx-c=0$. Then $ax_2-\frac{b}{2}=\sqrt{\frac{b^2}{4}+ac}$.  
Consider
$$
r^{i,j}(x)=\frac{-S^{i,j}(x)}{\sqrt{S^{i,i}S^{j,j}}}
$$
According to Silvester determinant identity one can write :
$$
S^{\{i,j\},\{i,j\}} \det S(x)=S^{i,i}S^{j,j}-[S^{i,j}(x)]^2
$$
Therefore for $x=x_1$ and $x=x_2$ one has
$$
S^{i,i}S^{j,j}-[S^{i,j}(x)]^2=0
$$
That is for $x=x_1$ and $x=x_2$ one has $r^{i,j}(x)=\pm 1$. The equation  $S^{i,j}(x)=ax-\frac{b}{2}$ implies that when $x$ is increasing from $x_1$ to $x_2$ then $r^{i,j}(x)$ is decreasing from $1$ to $-1$.
That is $r^{i,j}(x_2)=-1$, i.e.  $ax_2-\frac{b}{2}=\sqrt{S^{i,i}S^{j,j}}$. Therefore 
$$
\sqrt{S^{i,i}S^{j,j}}=\sqrt{\frac{b^2}{4}+ac}
$$
The Theorem is proved. 

Finally, the UMPU test for testing conditional independence of $X_i$ and $X_j$ can be written in the following form
\begin{equation}\label{Neyman_structure_final}
\varphi_{i,j}=\left\{\ 
\begin{array}{ll} 
0,& \displaystyle 2q-1 < r^{i,j} < 1-2q \\
1& \mbox{otherwise}
\end{array}\right.
\end{equation}
where $r^{i,j}$ is the sample partial correlation, and $q$ is the $\frac{\alpha}{2}$-quantile of beta distribution $Be(\frac{n-N}{2},\frac{n-N}{2})$,

\section{Concluding remarks}\label{conclusion}
In general optimality of tests for individual hypotheses testing does not imply optimality of multiple testing procedures. 
However if the losses from false decisions are supposed to be additive then it is possible to prove optimality from decision-theoretic point of view of some multiple testing procedures
\cite{bib10}, \cite{bib11}, \cite{bib12}. Application of this approach for Gaussian graphical model selection is a subject of further investigations.

\noindent
{\bf Acknowledgment:} The work was conducted at National Research University Higher School of Economics, Laboratory of Algorithms and Technologies for Network Analysis. 
Partly supported by  RFFI 14-01-00807.


\end{document}